\RequirePackage{fix-cm}

\documentclass[twoside]{article}
\usepackage{arxiv}

\usepackage{graphicx}

\usepackage{todonotes}
\usepackage{amsmath}
\usepackage{amsfonts}	
\usepackage{amssymb}
\usepackage[colorlinks,citecolor=blue,urlcolor=blue]{hyperref}
\usepackage{booktabs}

\DeclareTextFontCommand{\texttt}{\ttfamily\upshape}

\begin{document}

\title{The MIP Workshop 2023 Computational Competition on Reoptimization}

\author{%
	Suresh Bolusani$^*$\inst{1} \href{https://orcid.org/0000-0002-5735-3443}{\includegraphics[scale=0.1]{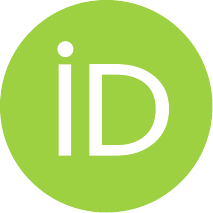}} \and
	Mathieu Besan{\c{c}}on\inst{1} \href{https://orcid.org/0000-0002-6284-3033}{\includegraphics[scale=0.1]{figures/orcid.pdf}} \and
	 Ambros Gleixner\inst{1,2} \href{https://orcid.org/0000-0003-0391-5903}{\includegraphics[scale=0.1]{figures/orcid.pdf}} \and
	Timo Berthold\inst{3} \href{https://orcid.org/0000-0002-6320-8154}{\includegraphics[scale=0.1]{figures/orcid.pdf}} \and
	Claudia D'Ambrosio\inst{4} \href{https://orcid.org/0000-0002-4040-0960}{\includegraphics[scale=0.1]{figures/orcid.pdf}} \and
	Gonzalo Mu{\~n}oz\inst{5} \href{https://orcid.org/0000-0002-9003-441X}{\includegraphics[scale=0.1]{figures/orcid.pdf}} \and
	Joseph Paat\inst{6} \href{https://orcid.org/0000-0003-2930-3155}{\includegraphics[scale=0.1]{figures/orcid.pdf}} \and
	Dimitri Thomopulos\inst{7} \href{https://orcid.org/0000-0003-0601-1790}{\includegraphics[scale=0.1]{figures/orcid.pdf}}
}

\institute{%
	Zuse Institute Berlin, Berlin, Germany \and
	HTW Berlin, Berlin, Germany \and
	TU Berlin \& FICO, Berlin, Germany \and
	LIX, CNRS, {\'E}cole Polytechnique, Institut Polytechnique de Paris, Palaiseau, France \and
	Universidad de O'Higgins, Rancagua, Chile \and
	University of British Columbia, Sauder School of Business, Vancouver, Canada \and
	Universit{\`a} di Pisa, DESTEC, Pisa, Italy
}

\email{bolusani@zib.de}

\date{29 November 2023}

\maketitle

\begin{abstract}
This paper describes the computational challenge developed for a computational competition held in 2023 for the $20^{\textrm{th}}$ anniversary of the Mixed Integer Programming Workshop.
The topic of this competition was reoptimization, also known as warm starting, of mixed integer linear optimization problems after slight changes to the input data for a common formulation.
The challenge was to accelerate the proof of optimality of the modified instances by leveraging the information from the solving processes of previously solved instances, all while creating high-quality primal solutions.
Specifically, we discuss the competition's format, the creation of public and hidden datasets, and the evaluation criteria.
Our goal is to establish a methodology for the generation of benchmark instances and an evaluation framework, along with benchmark datasets, to foster future research on reoptimization of mixed integer linear optimization problems.

\end{abstract}

\section{Introduction}\label{sec:intro}

The Mixed Integer Programming Workshop\footnote{\url{https://www.mixedinteger.org/\#mipworkshops}} is an annual single-track workshop highlighting the latest trends in mixed integer optimization and discrete optimization.
Since the computational development of optimization tools is a crucial component within the mixed integer linear optimization community, the workshop established a computational competition in 2022 to encourage and provide recognition to the development of novel practical techniques in the software for solving mixed integer linear optimization problems (MILPs).
The first edition of the competition focused on \emph{primal heuristics}, i.e., finding good quality primal solutions of general MILPs selected mainly from the MIPLIB 2017 benchmark library~\cite{miplib2017}.
This paper discusses the second edition of the competition held in 2023~\cite{MIPcc23Website}.
From an organizational point of view, this edition involved not only the development of the competition topic, its structure, and an evaluation framework but also the creation of new benchmarking instances.

Traditional benchmarks for MILPs often focus on the performance of optimizing a given instance from scratch. In many practically relevant settings, however, MILP solvers are used to repeatedly solve a series of similar instances of the same type and only slight modifications of the input data. This motivated the 2023 competition topic: the development of effective techniques for reoptimization, also known as warm starting, of MILPs in this setting.

In addition to the use case where practitioners solve the same MILP model with slightly perturbed input data, this setting also appears algorithmically:
\begin{itemize}
	\item Row generation algorithms, for instance based on generalized Benders' decompositions for single- and multilevel MILPs (e.g., bilevel optimization problems), solve MILP subproblems that vary only in the right-hand side vector across iterations, see, e.g.,~\cite{bolusani2022framework}.
	\item Column generation algorithms, for instance based on a Dantzig-Wolfe decomposition for MILPs, solve subproblems (pricing problems) that vary only in the objective vector across iterations, see, e.g.,~\cite{witzig2014reoptimization}.
	\item Scenario decomposition algorithms for stochastic MILPs solve subproblems that vary only in the scenario-dependent components both within and across iterations, see, e.g.,~\cite{Hassanzadeh2015}.
	  \item Primal heuristics, such as diving and neighborhood heuristics, may solve similar MILPs with varying input data both within and across calls to the heuristics, see, e.g.,~\cite{gamrath2019structure}.
\end{itemize}

Despite the broad applicability of reoptimization of MILPs, the research in this area is limited.
For example,~\cite{desrochers1988reoptimization} discusses a reoptimization algorithm for solving the shortest path problem with time windows in a dynamic programming setting,~\cite{ausiello2011complexity,schieber2018theory} develop frameworks for reoptimization of combinatorial optimization problems and derive specific algorithms for certain classes of such problems, and~\cite{ralphs2006duality,witzig2014reoptimization} discuss reoptimization techniques for general MILPs but only when specific components of the MILP vary from one instance to another.
Most of the existing literature in this area either deals with specific applications or needs to be more scalable to be applied in practice in the reoptimization settings mentioned above.
This further motivated the competition topic.

To evaluate reoptimization for the different settings mentioned above, the competition provided a set of MILP instance series, each of which comprised 50~related instances of the same size.
For each instance series, the type of change that may occur and the names of varying columns and/or rows as applicable for this type are known.

The competition participants were asked to provide a general solver to optimize a series of related MILPs in sequential order, thereby reusing information from the previous runs in the series.
The participants were free to build on any software that is available in source code, and that can be used freely for the evaluation of the competition.
The intention of this competition was not to perform offline training on different types of applications but to reuse information from current and previous solving processes to accelerate the solution of future instances.

The remainder of the paper has the following structure.
Section~\ref{sec:dataset} introduces the competition dataset, mentions the two GitHub repositories in which this dataset and its generation scripts are available, discusses the dataset creation process and the two metrics used for this purpose, and provides detailed explanations of all series of instances.
Section~\ref{sec:eval} discusses the competition's evaluation criteria and presents a novel scoring function developed to measure the computational performance of the submissions.
Finally, Section~\ref{sec:results} presents the competition results, concluding remarks, and future outlook.

\section{The Dataset}\label{sec:dataset}

The competition dataset consists of a set of MILP instance series of 50~related instances each.
Each instance series is based on an MILP taken from a specific application or benchmark library in the literature.
There are a total of seven public and five hidden instance series with constant constraint matrix, and three additional instance series with varying constraint matrix that were not part of the competition evaluation.

The instances in each series comply with the following specifications:
\begin{itemize}
	\item the number of constraints, and the number, order, and meaning of variables remain the same across the instances in a series, and
	\item some or all of the following input can vary:
	\begin{itemize}
		\item objective function coefficients,
		\item variable bounds,
		\item constraint sides, and
      \item coefficients of the constraint matrix.
	\end{itemize}
\end{itemize}

Due to the more challenging nature of constraint coefficient changes, the corresponding three series were not part of the official computational evaluation.
They were nonetheless included to see if some proposed approaches were applicable and efficient for these series.

Table~\ref{tab:dataset} summarizes all the instance series where
\begin{itemize}
	\item LO, UP, OBJ, LHS, RHS, and MAT, denote lower bound vector, upper bound vector, objective function vector, left-hand side vector, right-hand side vector, and the constraint matrix, respectively, and
	\item the last column states the time limit imposed to solve one instance in the series.
\end{itemize}

\begin{table}
\begin{center}
\caption{The Dataset}
\label{tab:dataset}
\newcommand{\xmark}{\checkmark}
\begin{tabular}{lccccccr}
\toprule
Instance series & \multicolumn{6}{c}{Varying component} & \quad Time limit \\
\cmidrule{2-7}
& LO & UP & OBJ & LHS & RHS & MAT \\
\midrule
\emph{Public}\\
\texttt{bnd\_series\_1} & & \xmark & & & & & 600\,s \\
\texttt{bnd\_series\_2} & \xmark & \xmark & & & & & 300\,s \\
\texttt{obj\_series\_1} & & & \xmark & & & & 400\,s \\
\texttt{obj\_series\_2} & & & \xmark & & & & 300\,s \\
\texttt{rhs\_series\_1} & & & & \xmark & \xmark & & 400\,s \\
\texttt{rhs\_series\_2} & & & & & \xmark & & 60\,s \\
\texttt{rhs\_obj\_series\_1} & & & \xmark & \xmark & \xmark & & 500\,s \\
\midrule
\emph{Hidden}\\
\texttt{bnd\_series\_3} & \xmark & \xmark & & & & & 600\,s \\
\texttt{obj\_series\_3} & & & \xmark & & & & 350\,s \\
\texttt{rhs\_series\_3} & & & & \xmark & \xmark & & 550\,s \\
\texttt{rhs\_series\_4} & & & & \xmark & \xmark & & 60\,s \\
\texttt{rhs\_obj\_series\_2} & & & \xmark & \xmark & \xmark & & 250\,s \\
\midrule
\emph{Out of competition}\\
\texttt{mat\_series\_1} & & & & & & \xmark & 300\,s \\
\texttt{mat\_rhs\_bnd\_series\_1} & \xmark & \xmark & & \xmark & \xmark & \xmark & 400\,s \\
\texttt{mat\_rhs\_bnd\_obj\_series\_1} & \xmark & \xmark & \xmark & \xmark & \xmark & \xmark & 180\,s \\
\bottomrule
\end{tabular}
\end{center}
\end{table}

To the best our knowledge, there is currently no dedicated reference benchmark for the reoptimization of MILPs.
We, therefore, created these series of instances based on various algorithms, applications, and existing MILP instances from the literature
and made them available at the GitHub repository~\cite{MIPcc23Website}.
Furthermore, the scripts used to generate these series are also available at another GitHub repository~\cite{MIPcc23DatasetGenScripts}, offering a template for generating more and longer series for future research beyond the context of the competition.

To build a single series, we generated numerous instances satisfying the series requirements, solved them to optimality from scratch, and selected 50~suitable instances for the series by applying the following two metrics.
\begin{enumerate}
\item The time taken to solve an instance to optimality.
\item The similarity between the varying components of the series, whenever applicable.
The similarity between two vectors $c$ and $\bar c$ is defined as the angle between the vectors~\cite{witzig2014reoptimization}.
\begin{equation*}
\text{Similarity} = \frac{ \langle c, \bar c \rangle}{\|c\| \|\bar c\|}
\end{equation*}
For example, let $c$ and $\bar c$ be the objective function vectors of two instances. If the similarity between these vectors is very high, then the effort required to solve these instances to optimality without reoptimization is likely to be similar.
Furthermore, most of the primal information generated in the solving process of the instance with the objective vector $c$ will be valid for the instance with the objective vector $\bar c$ (and vice versa).
\end{enumerate}
As mentioned earlier, the competition aims to encourage reoptimization techniques that can accelerate the solution of future instances in a given series compared to solving these instances from scratch.
Accordingly, we chose the time limit for an instance in a given series based on the individual solving times (from scratch) of the 50 instances in this series.
Then, the scoring function discussed in Section~\ref{sec:eval} penalized any techniques that could not solve an instance to optimality within this time limit.

We now discuss the details of each instance series.

\paragraph{\underline{\texttt{bnd\_series\_1}}}: This series is based on the instance \texttt{rococoC10-001000} from the MIPLIB 2017 benchmark library~\cite{miplib2010}. The instances were generated by perturbing the upper bounds of general integer variables selected via a discrete uniform distribution up to $\pm$100\% of the bound value.

\paragraph{\underline{\texttt{bnd\_series\_2}}}:
This series is based on the instance \texttt{csched007} from the MIPLIB 2017 benchmark library~\cite{miplib2010}. The instances were generated via random fixings of 15\% to 25\% of the binary variables selected via a discrete uniform distribution w.r.t.~the original instance.

\paragraph{\underline{\texttt{bnd\_series\_3}}}:
This series is also based on the instance \texttt{csched007} from the MIPLIB 2017 benchmark library~\cite{miplib2010}.
The instances were generated via random fixings of 5\% to 20\% of the binary variables selected via a discrete uniform distribution w.r.t.~the original instance.
Accordingly, these instances are relatively harder to solve as compared to the instances in the series \text{\texttt{bnd\_series\_2}} as indicated by the time limits in Table~\ref{tab:dataset}.

\paragraph{\underline{\texttt{obj\_series\_1}}}:
This series is based on the stochastic multiple binary knapsack problem and the associated instance set introduced in~\cite{angulo2016improving}.
The problem has the formulation
\begin{align*}
	&& \min_{x,z} &\; c^\top x + d^\top z + Q(x)\\
        && \text{s.t.} & \; Ax + Cz \geq b,\\
        &&& \; x \in \{0, 1\}^n,\; z \in \{0, 1\}^n,
\end{align*}
where $Q(x) = \sum\limits_{\omega \in \Omega} p_\omega Q_{\omega}(x)$, with
\begin{align*}
	&& Q_{\omega}(x) := \min_{y} &\; q_{\omega}^\top y\\
        && \text{s.t.} & \; Wy \geq h - Tx,\\
        &&& \; y \in \{0, 1\}^n,
\end{align*}
where $\omega \in \Omega$ denotes a scenario, $p_\omega$ denotes probability of the scenario $\omega$, and the second-stage objective vector $q_{\omega}$ is random, following a discrete distribution with finitely many scenarios.

We adapted the given dataset and generated instances by considering one scenario at a time. This resulted in a series of instances with one-third of the objective vector (corresponding to $y$ variables) varying across instances.

\paragraph{\underline{\texttt{obj\_series\_2}}}: This series is based on the instance \texttt{ci-s4} from the MIPLIB 2017 benchmark library~\cite{miplib2017}. The instances were generated via random perturbations and random rotations of the objective vector.

\paragraph{\underline{\texttt{obj\_series\_3}}}: This series is based on the UCI Machine Learning repository dataset \texttt{magic}~\cite{misc_magic_gamma_telescope_159}. The instances are subproblems of a column generation algorithm for improving decision trees~\cite{FIRAT2020104866}. The final set of instances were generated based on a submission that was received in response to a public call for additional datasets.

\paragraph{\underline{\texttt{rhs\_series\_1}}}: This series is based on the stochastic server location problem and the associated dataset proposed in~\cite{ntaimo2010disjunctive}. The problem has the formulation
\begin{align*}
	&& \min_{x} &\; \sum\limits_{i \in {\mathcal{I}}} c_i x_i - Q(x)\\
        && \text{s.t.} & \; \sum\limits_{i \in {\mathcal{I}}} x_i \leq u,\\
        &&& \; x_i \in \{0, 1\}\; \forall i \in {\mathcal{I}},
\end{align*}
where $Q(x) = \sum\limits_{\omega \in \Omega} p_\omega Q_{\omega}(x)$, with
\begin{align*}
	&& Q_{\omega}(x) := \min_{y} &\; \sum\limits_{i \in {\mathcal{I}}} \sum\limits_{j \in {\mathcal{J}}} (-q_{ij} y_{ij}) + \sum\limits_{i \in {\mathcal{I}}} g_{i0} y_{i0}\\
        && \text{s.t.} & \; \sum\limits_{j \in {\mathcal{J}}} d_{ij} y_{ij} - y_{i0} \leq \tau_i x_i,\; \forall i \in {\mathcal{I}},\\
        &&& \sum\limits_{i \in {\mathcal{I}}} y_{ij} = r_{j\omega},\; \forall j \in {\mathcal{J}},\\
        &&& \; y_{ij} \in \{0, 1\}, \; \forall i \in {\mathcal{I}}, \forall j \in {\mathcal{J}},\; y_{i0} \geq 0,\; \forall i \in {\mathcal{I}},
\end{align*}
where $\omega \in \Omega$ denotes a scenario and $p_\omega$ denotes probability of the scenario $\omega$.

We adapted the given dataset and generated instances by considering 25~scenarios at a time. This resulted in a series of instances with only the right-hand side vector of equality constraints varying across instances.

\paragraph{\underline{\texttt{rhs\_series\_2}}}: This series is based on a synthetic MILP and the associated dataset proposed in~\cite{jimenez2022warm}. The problem has the following formulation.
\begin{align*}
	&& \min_{x,y} &\; c^\top x\\
        && \text{s.t.} & \; Ax \leq b,\\
        &&& \; l^\top y \leq x \leq u^\top y, \\
        &&& \; x \in \mathbb{R}^n,\; y \in \{0, 1\}^n
\end{align*}
We adapted the given dataset and generated instances by taking a convex combination of two different RHS vectors.

\paragraph{\underline{\texttt{rhs\_series\_3}}}: This series is based on the instance \texttt{glass4} from the MIPLIB 2017 benchmark library~\cite{miplib2003}. The instances were generated by perturbing the non-negative LHS and RHS vector components selected via a discrete uniform distribution up to $\pm 70\%$ of their values.

\paragraph{\underline{\texttt{rhs\_series\_4}}}: This series is also based on the synthetic MILP and the associated dataset proposed in~\cite{jimenez2022warm}. We adapted the given dataset and generated instances by taking a convex combination of two different RHS vectors (different than the ones used for generating \texttt{rhs\_series\_2}).

\paragraph{\underline{\texttt{mat\_series\_1}}}: This series is based on the optimal vaccine allocation problem and the associated dataset proposed in~\cite{tanner2010iis}. The problem formulation is
\begin{align*}
	&& \min_{x} &\; \sum\limits_{i \in {\mathcal{I}}} \sum\limits_{j \in {\mathcal{J}}} \sum\limits_{k \in {\mathcal{K}}} c_{ik} x_{ij}\\
        && \text{s.t.} & \; \sum\limits_{j \in {\mathcal{J}}} x_{ij} = 1,\; \forall i \in {\mathcal{I}},\\
        &&& \; \sum\limits_{i \in {\mathcal{I}}}\sum\limits_{j \in {\mathcal{J}}} a_{ij\omega} x_{ij} \leq 1 + M z_{\omega},\; \forall \omega \in {\Omega},\\
        &&& \; \sum\limits_{\omega \in {\Omega}} p_{\omega} z_{\omega} \leq 1 - \alpha,\\
        &&& \; x_{ij} \in [0, 1], \; \forall i \in {\mathcal{I}}, \forall j \in {\mathcal{J}},\; z_{\omega} \in \{0, 1\},\; \forall \omega \in \Omega,
\end{align*}
where $\omega \in \Omega$ denotes a scenario, $p_\omega$ denotes probability of the scenario $\omega$, and $M$ denotes a big-M parameter.

We adapted the given dataset and generated instances by considering 500~scenarios at a time. This resulted in a series of instances with the constraint matrix corresponding to inequality constraints varying across instances.

\paragraph{\underline{\texttt{rhs\_obj\_series\_1}}}: This series is based on the hydro unit commitment (HUC) problem modeled as an MILP. Considering a fixed hydro valley, the input data potentially changing is restricted to the electricity prices, the inflows, and the initial and target water volume in the reservoirs. These appear only in the objective function or constraint sides. Thus, reoptimizing this problem is practically interesting because the great majority of the input data remains unchanged. Moreover, utility companies often solve the HUC problem as a subproblem of a decomposition method. Consequently, for converging to the optimal solution of the whole unit commitment problem, a HUC has to be solved at each iteration.
Detailed mathematical formulation is available in~\cite{thomopulos2023generating}.

\paragraph{\underline{\texttt{rhs\_obj\_series\_2}}}: This series is based on a hydroelectric valley (Ain River) industrial use case in France. Six dams and their different turbines are modeled for the next four days with an hourly time step. The differences across instances are electricity prices (in the objective function) and the varying flows of the different affluents (in the RHS vector) of the river Ain, which were collected from~\cite{andreassian_2021,gme_2022}. This series is based on the instances that were received in response to a public call for additional datasets. The final instances were generated by perturbing the LHS, RHS, and objective function vector components, selected via a discrete uniform distribution up to $\pm$20\% of their values.

\paragraph{\underline{\texttt{mat\_rhs\_bnd\_series\_1}}}: This series of instance are based on the MILP formulation of the multilevel supply chain of a fictitious cell phone company. Detailed description is available at~\cite{cellphonedummy}.

\paragraph{\underline{\texttt{mat\_rhs\_bnd\_obj\_series\_1}}}: This series of instances is also based on the hydro unit commitment problem similar to the series \texttt{rhs\_obj\_series\_1}, but every data component of the given MILP can vary here.

\section{Evaluation Criteria}\label{sec:eval}
This section discusses the criteria used to evaluate the proposed approaches for the competition. Participants were asked to provide a written report and code base corresponding to their submission. Then, the following two criteria were used for the final evaluation.

\begin{enumerate}
	\item \textbf{Novelty and scope}: innovativeness of the approach and its general applicability w.r.t.~the varying components and magnitude of their variation.
	\item \textbf{Computational excellence}: ranking of the approach in terms of the performance score defined later in this section.
\end{enumerate}
\noindent
The written report had to include the following:
\begin{itemize}
	\item description of the methods developed and implemented, including any citations to the literature and software used,
	\item computational results on the public instance series with constant constraint matrix,
	\item analysis with at least $\texttt{reltime}$, $\texttt{gap}$, and $\texttt{nofeas}$ scores (defined later in the section) averaged over all 50 instances and additionally over the five batches of instances 1 to 10, 11 to 20, 21 to 30, 31 to 40, and 41 to 50, and
	\item any further analysis including the applicability of the approach to the instance series with varying constraint matrix.
\end{itemize}

Participants were allowed to use any existing software available in source code and freely usable for the evaluation.
A submission was required to solve the instances of one series sequentially in the order specified by the input files.
It was not permitted to parse and analyze instances in one series ahead of time, i.e., while solving the $i^\textrm{th}$ instance in the series, a submission may use only information from the first $i-1$ instances.
A submission was not allowed to modify the solution for an instance after moving on to solve the next instance. A submission must run sequentially (1 CPU thread), use no more than 16 GB of RAM, and respect the total time limit for each instance series.
Violations of the time limit for a single instance are penalized in the performance score.

Evaluating and comparing the performance of submissions necessitated the development of an appropriate scoring metric. Our goal was to create a comprehensive scoring system that takes into account various potential situations, ranked from the most favorable to the least favorable:

\begin{itemize}
\item the instance is solved to optimality within a certain runtime, shorter runtimes are preferred;
\item the instance times out but provides an incumbent solution, smaller gaps are preferred;
\item the instance times out without yielding any feasible solution, which is generally regarded as undesirable.
\end{itemize}

While some existing measures, such as~\cite{berthold2021confined}, work across these different situations, we also sought to reflect the aforementioned hierarchy of favorability in our scoring. It was crucial to strike a balance; we wanted to penalize, to a certain extent, failure to achieve optimality or to find a solution, but we also aimed to avoid unduly restricting approaches that engage in an early exploration phase, which might intentionally deteriorate performance for the first few instances of a series.

To achieve this balance, we established a scoring range:
\begin{itemize}
\item the instances solved to optimality within the time limit were assigned scores between 0 and 1,
\item the instances that timed out with an incumbent solution received scores between 1 and 2, and
\item the instances that timed out without any feasible solution were given a score of 3.
\end{itemize}

Specifically, let $s = 1, 2, \dots, S$ denote the index of the instance series, each consisting of $i = 1, 2, \dots, 50$ MILP instances. Then, the performance on a single instance $(s, i)$ is measured via the scoring function
\begin{equation*}
f(s, i) = \texttt{reltime} + \texttt{gap} + \texttt{nofeas},
\end{equation*}
where
\begin{equation*}
\texttt{reltime} = 
\begin{cases}
\frac{\texttt{time spent}}{\texttt{time limit}}, & \text{if the instance is solved to optimality},\\
\text{max}\left\{1, \frac{\texttt{time spent}}{\texttt{time limit}}\right\}, & \text{otherwise},
\end{cases}
\end{equation*}
\begin{equation*}
\texttt{gap} = 
\begin{cases}
0, & \text{if primal bound (\texttt{pb}) and dual bound (\texttt{db}) are 0},\\
1, & \text{if \texttt{pb} or \texttt{db} are infinite, or \texttt{pb} $\times$ \texttt{db} $< 0$},\\
\frac{|\texttt{pb} - \texttt{db}|}{\text{max}\left(|\texttt{pb}|, |\texttt{db}|\right)}, & \text{otherwise},
\end{cases}
\end{equation*}
and
\begin{equation*}
\texttt{nofeas} = 
\begin{cases}
1, & \text{if no feasible solution is returned},\\
0, & \text{otherwise}.
\end{cases}
\end{equation*}
Smaller scores $f(s, i)$ are better.
In our mind, this scoring framework allows for nuanced evaluation, rewarding efficiency while maintaining fairness across different methodologies. Note that the transition between solving an instance and timing out with a small gap is ``smooth'' in the sense that the score approaches 1 from below as the solve time gets close to the time limit and 1 from above as the final gap approaches 0 for unsolved instances. In contrast, the term \texttt{nofeas} adds a fixed penalty term for not producing a primal solution.

For technical reasons, we also had to address the situations where submissions exceed the time limit (usually only slightly) or stop a solution process prematurely before reaching the time limit.
If the submission exceeds the time limit, \texttt{reltime} will be larger than 1. Submissions that stop before the time limit without reaching zero gap will receive \texttt{reltime} = 1. These considerations were mainly taken to close loopholes for the competition, such as intentionally ignoring a time limit to gather more information or aborting unpromising runs to keep the score low.

Then, for each instance $(s, i)$, all participants were ranked according to their score $f(s, i)$, resulting in each team receiving a rank $r(s, i)$, where smaller is better. Teams with the same score received the same rank. For instances for which the primal solution was not feasible or the dual bound was not valid, a team received two times the worst possible rank independently of their $f(s, i)$ value.

Following the motivation of the challenge to reward methods that use information from
previous solving processes in order to gain performance, we assigned a gradually
increasing weight to later instances in each series, i.e., we computed the final score as
\begin{equation*}
C = \sum\limits_{s = 1}^{S} \sum\limits_{i = 1}^{50} (1 + 0.1 i)\;r(s, i).
\end{equation*}
The lower this final score, the better.

\section{Results and Outlook}\label{sec:results}

The jury selected one winner and awarded one honorable mention from a total of twelve registrations and six final submissions containing both a written report and implementation code.

The winning submission \emph{Progressively Strengthening and Tuning MIP Solvers for Reoptimization} by Krunal Patel described in detail in \cite{Patel2023} convinced the jury not only through its computational excellence that was displayed by the top ranked performance on almost all public and hidden data sets, but also through its broad applicability and attention to algorithmic details.
Building on top of an existing, open-source LP-based branch-and-bound solver \textsc{SCIP}~\cite{BestuzhevaEtal2021ZR}, the approach distinguished itself by targeting multiple aspects of the solving process in combination,
reusing primal information and pseudo costs and improving parameter configurations online despite the limited number of observations.

An honorable mention was awarded to the submission \emph{Influence branching for learning to solve mixed integer programs online} by the team of Paul Strang, Zacharie Ales, Come Bissuel, Olivier Juan, Safia Kedad-Sidhoum, and Emmanuel Rachelson.
The jury was impressed by the successful adaptation of influence branching~\cite{EtheveEtAl2020} to the reoptimization setting of the competition.
The submission employed online hyperparameter tuning of different influence models via multi-armed bandit selection and consistently performed well on both public and hidden datasets.

With the created instance series available on the repository~\cite{MIPcc23Website} and presented in Section~\ref{sec:dataset},
and their extensibility based on the generation scripts in~\cite{MIPcc23DatasetGenScripts},
we hope to offer the research community a set of benchmarks to foster and evaluate future research efforts towards reoptimization.
We aim to continuously develop~\cite{MIPcc23DatasetGenScripts} as a benchmarking library in the long run for the research on the reoptimization of MILPs.
Accordingly, we are expanding the current dataset by generating additional series of instances and making them available to the research community via this repository.

\begin{acknowledgements}
  We wish to thank
  Daniel Bienstock for providing us with computational infrastructure in order to evaluate the submissions,
  Dimitri Kniasew from SAP for submitting one hidden series of supply chain network planning instances,
  Felipe Serrano for discussing initial suggestions for the topic of the competition, and
  Domenico Salvagnin for participating in the final evaluation process.
\end{acknowledgements}

\section*{Conflict of Interest}
The authors declare that they have no conflict of interest.

\bibliographystyle{spmpsci}      
\bibliography{references}

\end{document}